\documentclass[10pt]{article}

\usepackage{amsfonts}
\usepackage{amsmath}

\newcommand{\vp}{\varphi}

\begin{document}

\begin{center}
 {\large \bf Direct and inverse scattering on noncompact star-type quantum graph with Bessel singularity
 } \\[0.2cm]
 {\bf Mikhail Ignatyev} \\[0.2cm]
\end{center}

{\bf Abstract.}
In this paper we study the noncompact star-type graph with perturbed radial Schrodinger equation on each ray and the matching conditions of some special form at the vertex. The results include the uniqueness theorem and constructive procedure for solution of the inverse scattering problem.
\\
[0.1cm]



\section{Introduction}

Transport, spectral and scattering problems
for differential operators on graphs appear frequently in mathematics,
natural sciences and engineering {\cite{FaP}, \cite{Ex}, \cite{KoS},
 \cite{PoB}, {\cite{Kuch}, \cite{AKu}} and are in the focus of intensive
investigations {\cite{Ger}},{\cite{BrW}},  {\cite{Bel1}}, {\cite{Yur1}}, {\cite{Yur2}}, {\cite{Troo}}, {\cite{KuS}},{\cite{Kur}}, {\cite{ISI}}, {\cite{Tro1}}.

In this paper we study the noncompact star-type graph with perturbed radial Schrodinger equation  on each ray. Due to its physical importance this equation has obtained much attention and we refer for example to {\cite{LD}}, {\cite{Coz76}}, {\cite{Coz83}}, {\cite{KoTeS}}, {\cite{Hr12}}, {\cite{KoTe}}  and the references therein. Our studies dealing with the equation on a metric graph can be treated as some continuation of the works {\cite{Yur3}}, {\cite{Fed}} (see also the references therein), where the important class of the Sturm--Liouville operators having Bessel-type singularity inside the interval was investigated. Matching conditions at the vertex arise as a generalization of matching conditions introduced in the works above under some additional conditions on behavior of the potential near the singularity. We call them the {\it matching conditions of analytical type }. One can notice that even for Sturm--Liouville expressions with real-valued potentials on finite interval such matching conditions can generate nonself-adjoint operators. For this reason we do not restrict our considerations to the case of real-valued potentials and real coefficients in matching conditions.

Our considerations follow in general the scheme recently developed for the investigation of inverse spectral problems on quantum graphs. The key  point is the solution of "partial inverse problem" consisting in recovering the potential on one ray from the certain part of scattering data. Although in technical sense such problem looks similar to the classical inverse scattering problem for radial Schrodinger equation on the half-line, our approach differs from the approach used in the previous works {\cite{Coz83}}. Instead of classical Marchenko method we use a certain version of contour integral method based in part on the approach presented in {\cite{Bea}} and in part on the ideas of spectral mappings method {\cite{YuB}}. This allows us to provide the constructive procedure for solution the inverse problem based a certain linear equation with less restrictive requirements on the behavior of the potential at infinity; moreover, contrary to {\cite{Coz83}}, we are able to prove the solvability of this equation. More exactly, under some conditions of regularity and "genericity" (which are similar to "genericity" conditions introduced in {\cite{Bea}} and assume the absence of spectral singularities and some kind of regularity for small values of spectral parameter) our method works provided the potential satisfies at infinity only Marchenko decay condition.

Main results of the paper are contained in Theorem 5.1 and include the uniqueness result and constructive procedure for solution of the inverse scattering problem.

\section{Direct scattering on the semi-axis. Auxiliary facts and notations.}

Consider the equation:
$$
-y''+\left(\frac{\nu_0}{x^2}+q(x)\right)y=\lambda y=\rho^2 y, \ x>0 \eqno(2.1)
$$
with complex $\nu_0$ and complex-valued $q$ that will be called as a {\it potential}. Let $\nu_0=\nu^2-1/4$. For definiteness we assume that $\mbox{Re}\nu>1/2$, $\nu\notin\mathbb{N}$ and $q$ satisfies the conditions:
$$
\int\limits_0^1 \left|x^{1-2\nu} q(x)\right|dx +
\int\limits_1^\infty \left|xq(x)\right|dx<\infty. \eqno(2.2)
$$
Here and below we assume $z^\mu:=\exp(\mu\log|z|+i\mu\arg z)$, $\arg z\in(-\pi,\pi]$.
Together with (2.1) we consider the unperturbed equation:
$$
-y''+\frac{\nu_0}{x^2}y=\lambda y=\rho^2 y, \ x>0 \eqno(2.3)
$$

Denote $\mu_1=1/2-\nu$, $\mu_2=1/2+\nu$. Let $C_j(x,\lambda)$, $j=1,2$ be the solutions of unperturbed equation (2.3) defined as:
$$
C_j(x,\lambda)=x^{\mu_j}\sum\limits_{k=0}^\infty c_{jk}\lambda^k x^{2k}, \ c_{10}c_{20}=(2\nu)^{-1},
$$
$$
c_{jk}=(-1)^k c_{j0}\left(\prod\limits_{s=1}^k((2s+\mu_j)(2s+\mu_j-1)-\nu_0)\right)^{-1}.
$$
$C_j(x,\lambda)$ are both entire with respect to $\lambda$.

Define $S_j(x,\lambda)$, $j=1,2$ as the solutions of (2.1) satisfying the integral equations:
$$
S_j(x,\lambda)=C_j(x,\lambda)-\int\limits_0^x g(x,t,\lambda) q(t)S_j(t,\lambda) dt,
$$
where $g(x,t,\lambda)$ is a Green function: $g(x,t,\lambda)=C_1(x,\lambda)C_2(t,\lambda)-C_2(x,\lambda)C_1(t,\lambda)$. Note that the equations are solvable for both $j=1,2$ provided the condition (2.2).
 We also consider the Jost solution $f(x,\rho)$ for (2.1) normalized with the asymptotics $f^{(\xi)}(x,\rho)=\mbox{e}^{i\rho x}((i\rho)^{\xi}+o(1))$, $x\to\infty$, $\xi=0,1$. Such solution exists (and unique) under the condition $xq(x)\in L_1(1,\infty)$, analytic with respect to $\rho\in\Omega_+:=\{\rho:\mbox{Im}\rho>0\}$ and $\rho^{-\mu_1}f^{(\xi)}(x,\rho)$ is continuous in $\rho\in\overline\Omega_+$.

The following expansion plays an important role in our further considerations:
$$f(x,\rho)=b_1(\rho)S_1(x,\lambda)+b_2(\rho)S_2(x,\lambda),$$
where coefficients $b_j(\rho)$, $j=1,2$ are called the {\it Stokes multipliers}. The ratio $b_2(\rho)/b_1(\rho)$ can be easily shown to coincide with the Weyl function $m(\lambda)$, which determines uniquely the potential $q$. This result and constructive procedure for solving the corresponding inverse problem one can find, for instance, in {\cite{Yur4}}.

The following properties of $b_j(\rho)$ arise from the properties of Jost solution $f(x,\rho)$ and the results of {\cite{Yur4}}.

\medskip
{\bf Lemma 2.1.}  For $\rho\to\infty$ the following asymptotics holds:
$$
b_j(\rho)=\rho^{\mu_j}\left(b_j^\infty+O(\rho^{-1})\right),
$$
where the constants $b_j^\infty\neq 0$ do not depend on $q$.

For $\rho\to 0$ the following asymptotics holds:
$$
b_j(\rho)=\rho^{\mu_1}\left(b_j^0+o(1)\right)
$$
with some constants $b_j^0$ (depending, in general, on $q$).

\section{Direct scattering on the graph}

We consider the noncompact star-type metric graph $\Gamma$ consisting of finite number of rays $\{\mathcal{R}_k\}_{k=1}^p$ emanating from their common initial vertex. A function $y$ on the ray $\mathcal{R}_k$ we consider as a function of local parameter $x\in[0,\infty)$, where $x=0$ corresponds to the vertex. A function $y$ on $\Gamma$ we consider as a set of functions $\{y_k\}_{k=1}^p$ (where $y_k=y\left.\right|_{\mathcal{R}_k}$ is considered as a function on $[0,\infty)$ as it was specified above).

On each ray $\mathcal{R}_k$, $k=\overline{1,p}$ we consider the differential equation
$$
-y''+\left(\frac{\nu_{k0}}{x^2}+q_k(x)\right)y=\lambda y=\rho^2 y. \eqno(3.1)
$$
We follow the terminology of previous section and use the same notations adding the number $k$ of ray as a first index. In particular, we assume $\nu_{k0}=\nu_k^2-1/4$, $\mbox{Re}\nu_k>1/2$, $\nu_k\notin\mathbb{N}$ and $q_k(x)$ to satisfy (2.2) (with $\nu=\nu_k$). The complex numbers $\nu_{k0}$ may be different for different rays but the following restriction is assumed throughout the paper.

\medskip
{\bf Condition $\nu$.} If $\nu_j\neq\nu_k$ then $\mbox{Re}\nu_j\neq\mbox{Re}\nu_k$.

\medskip
Let some function $y$ satisfy (3.1). Then the Wronskians $\langle S_{k1},y\rangle$ and $\langle y,S_{k2}\rangle$ do not depend on $x$ and we can introduce the following linear forms:
$$
U_{k1}(y):=\sigma_k\langle y,S_{k2}\rangle,
$$
$$
U_{k2}(y):=\sigma_{k1}\langle y,S_{k2}\rangle+\sigma_{k2}\langle S_{k1},y\rangle.
$$
Everywhere in the sequel we assume $\sigma_k\neq 0, \sigma_{k2}\neq 0$, $k=\overline{1,p}$.

Now let $y=\{y_k\}_{k=1}^p$ be a function on $\Gamma$ such that each of $y_k$ satisfies (3.1). We define the {\it matching conditions} at the vertex of $\Gamma$ as follows:
$$
U_{j1}(y_{j})=U_{k1}(y_{k}),  j\neq k, \ \sum\limits_{j=1}^p U_{j2}(y_{j})=0. \eqno(3.2)
$$
Matching conditions (3.2) are generalization of matching conditions introduced in {\cite{Yur3}} for Sturm-Liouville operators having singularities inside the interval. They also can be considered as some generalization of standard (i.e. continuity together with Kirchhoff condition) matching conditions for quantum graphs.

With each ray $\mathcal{R}_k$ we associate the{\it Weyl-type solution} $\psi_k(\rho)$ that we determine as a function on $\Gamma$ $\psi_k(\rho)=\left\{\psi_{kj}(x,\rho)\right\}_{j=1}^p$ such that:
\begin{itemize}
\item each of $\psi_{kj}$, $j=\overline{1,p}$ solves the differential equation (3.1);
\item $\psi_{kj}(x,\rho)=O\left(\mbox{e}^{i\rho x}\right)$ as $x\to\infty$ for $j\neq k$;
\item $\psi_{kk}(x,\rho)=\mbox{e}^{-i\rho x}(1+o(1))$ as $x\to\infty$;
\item $\psi_k(\rho)$ satisfies the matching conditions (3.2).
\end{itemize}

Let $k$ be arbitrary fixed. We look for $\psi_k$ in the following form:
$$
\psi_{kj}(x,\rho)=\gamma_{kj}(\rho)f_j(x,\rho), j\neq k, \ \psi_{kk}(x,\rho)=\gamma_{kk}(\rho)f_k(x,\rho)+\frac{2i\rho}{b_{k1}(\rho)} S_{k2}(x,\lambda). \eqno(3.3)
$$
It is clear that any function of the form (3.3) with arbitrary coefficients $\gamma_{kj}$ satisfies all the conditions determining the Weyl-type solution except matching conditions (3.2). Substituting (3.3) into (3.2) and taking into account that $\langle f_j, S_{j2}\rangle=b_{j1}$, $\langle  S_{j1},f_j\rangle=b_{j2}$ we arrive at the following linear algebraic system:
$$
\sigma_j b_{j1}\gamma_{kj}+\beta_k=0, \ j=\overline{1,p}, \ \sum\limits_{j=1}^p(\sigma_{j1} b_{j1}+\sigma_{j2} b_{j2})\gamma_{kj}=-\sigma_{k2} \delta_k, \ \delta_k=\frac{2 i\rho}{b_{k1}} \eqno(3.4)
$$
with respect to the values $\gamma_{kj}, j=\overline{1,p}$, $\beta_k$. Cramer rule being applied to (3.4) together with representation (3.3) shows that  $\psi_k(\rho)$ exists (and is unique) for all $\rho\in \overline \Omega_+$ outside the zeros of the function $b_{k1}(\rho)\Delta(\rho)$, where $$\Delta(\rho)=\sum\limits_{s=1}^p(\sigma_{s1} b_{s1}(\rho)+\sigma_{s2} b_{s2}(\rho))\prod\limits_{j\neq s}\sigma_j b_{j1}(\rho) \eqno(3.5)$$
is a determinant of the system (3.4). In the sequel we assume that the following "genericity" condition holds.

\medskip
{\bf Condition $G$.} All the functions $b_{k1}(\rho)\Delta(\rho)$, $k=\overline{1,p}$ have no real zeros (with possible exception of $\rho=0$) and no multiple zeros.

Let $\{1,...,p\}=\bigcup\limits_{\xi=1}^m I_\xi$, where the pairwise disjoint subsets $I_\xi$ are such that for any $\xi$, $j,k\in I_\xi$ one has $\nu_j=\nu_k=\tau_\xi$ and $\mbox{Re}\tau_1>\dots\mbox{Re}\tau_m$. The following properties of characteristic function $\Delta(\rho)$ are the direct sequences of representation (3.5) and lemma 2.1.

\medskip
{\bf Lemma 3.1.} For $\rho\to\infty$ the following asymptotics holds:
$$\Delta(\rho)=\sum\limits_{\xi=1}^m\rho^{\mu_1+2\tau_\xi}\left(a^\infty_\xi+O(\rho^{-1})\right),$$
where $$a^\infty_\xi=\sum\limits_{s\in I_\xi}\sigma_{s2} b^\infty_{s2}\prod\limits_{j\neq s}\sigma_j b^\infty_{j1}, \ \mu_1:=\sum\limits_{j=1}^p \mu_{j1}.$$

For $\rho\to 0$ one has: $$\Delta(\rho)=\rho^{\mu_1}(d^0+o(1)), \ d^0=\sum\limits_{s=1}^p(\sigma_{s1} b^0_{s1}+\sigma_{s2} b^0_{s2})\prod\limits_{j\neq s}\sigma_j b^0_{j1}.$$

\medskip
It is important to notice that the coefficients $a^\infty_\xi$ depend on the numbers $\nu_k$, $k=\overline{1,p}$ and coefficients $\sigma_k, \sigma_{k1}, \sigma_{k2}$ in matching conditions (3.2) and do not depend on the potentials $q_k$.
In what follows we assume that the following conditions of "regularity" and "genericity at 0" hold.

\medskip
{\bf Condition $R_\infty$.} $a^\infty_1\neq 0$.

\medskip
{\bf Condition $R_0$.} $b^0_{j1}\neq 0$, $j=\overline{1,p}$, $d^0\neq 0$.

\medskip
Our further considerations deal with $\psi_{kk}(x,\rho)$. Using (3.3) and the
representations:
$$
\gamma_{kk}(\rho)=\delta_k(\rho)\cdot\frac{\Delta_k(\rho)}{\Delta(\rho)},\ \delta_k(\rho)=\frac{2i\rho}{b_{k1}(\rho)},\ \Delta_k(\rho)=(\sigma_k)^2\prod\limits_{j\neq k}\sigma_j b_{j1}(\rho), \eqno(3.6)
$$
that can be obtained from (3.4) using the Cramer rule we deduce the following assertion.

\medskip
{\bf Lemma 3.2.} Under the Condition $R_0$ one has $$\gamma_{kk}(\rho)=O\left(\rho^{1-2\mu_{k1}}\right)$$ as $\rho\to 0$. Moreover, $$\psi^{(\xi-1)}_{kk}(x,\rho)=O\left(\rho^{\mu_{k2}}\right), \ \xi=0,1$$ as $\rho\to 0$ for any fixed $x>0$.

\medskip
It follows from lemmas 2.1, 3.1 that under the conditions $G$, $R_0$ and $R_\infty$ the function $b_{k1}(\rho)\Delta(\rho)$ has at most finite set of simple zeros, $\psi_{kk}(x,\rho)$ has at these points either a simple pole or a removable singularity. Denote the set of poles of $\psi_{kk}(x,\rho)$ in $\Omega_+$ as $Z^+_k$. It follows from representation (3.3) that for any $\rho_0\in Z^+_k$ one has:
$$
{\rm res}_{\rho=\rho_0}\psi_{kk}(x,\rho)=\alpha_k(\rho_0) f_k(x,\rho_0)=\alpha_k(\rho_0) \exp(i\rho_0 x)(1+o(1)), x\to\infty. \eqno(3.7)
$$
with some constant $\alpha_k(\rho_0)$. Indeed, (3.7) is obvious if $\rho_0$ is a zero of $\Delta(\rho)$; if $\rho_0$ is a zero of $b_{k1}(\rho)$ we take into account that $S_{k2}(x,\rho_0)$ is proportional to $f_k(x,\rho_0)$ and arrive again to (3.7) (let us recall that under the condition $G$ $\Delta(\rho)$ and $b_{k1}(\rho)$ have no common zeros).

Further, for real nonzero $\rho$ usual arguments yield the following representation:
$$
\psi_{kk}(x,\rho)=f_k(x,-\rho)+r_k(\rho)f_k(x,\rho)=\exp(-i\rho x)+r_k(\rho)\exp(i\rho x)+o(1), x\to\infty. \eqno(3.8)
$$
We call $r_k$ the {\it reflection coefficient associated with the ray $\mathcal{R}_k$}, the set $J_k:=\{r_k, Z^+_k, \alpha_k(\rho), \rho\in Z^+_k\}$ is called the {\it scattering data associated with the ray $\mathcal{R}_k$}, the set $J=\{J_k\}_{k=1}^{p-1}$ is called the {\it scattering data for $\Gamma$}.

\section{Inverse scattering. Partial problem.}

This section is devoted to the following partial inverse problem.

\medskip
{\bf Problem $IP(k)$.} Given $J_k$, $k\in\{1,\dots,p\}$ find $q_k(x), x>0$.

\medskip
Let $L$ denotes the problem on $\Gamma$ consisting of differential equations (3.1) on each ray $\mathcal{R}_j$, $j=\overline{1,p}$ and matching conditions (3.2). Together with $L$ we consider the problem $\tilde L$ with equations of the same form (3.1) and matching conditions (3.2) but different potentials $\tilde q_j$, $j=\overline{1,p}$. We assume that coefficients in (3.2) are the same for $L$ and $\tilde L$ and $\tilde \nu_{k0}=\nu_{k0}$. Moreover, we assume that the conditions $G$, $R_0$ and $R_\infty$ hold for both $L$ and $\tilde L$. We agree that if $\eta$ denotes some object related to $L$ then $\tilde \eta$ denotes an analogous object related to $\tilde L$ and $\hat\eta:=\eta-\tilde\eta$.

In what follows we shall also use the following notations.  If $A$ denotes some matrix then $A_j$ denotes its $j$-th row. If $f$ is some function holomorphic in deleted neighborhood of $\rho_0$ then $f_{\langle m\rangle}(\rho_0)$ denote the coefficients in the Laurent expansion:
$$
f(\rho)=\sum\limits_{m=-\infty}^\infty (\rho-\rho_0)^m f_{\langle m\rangle}(\rho_0).
$$
If $f(\rho)$ is some function meromorphic outside the real axis then for real $\rho$ we denote $f_\pm(\rho):=\lim\limits_{\varepsilon\to +0}f(\rho\pm i\varepsilon)$ .

From this point in this section we assume that $k$ is arbitrary fixed. The key role in our further considerations is played by the  {\it spectral mappings matrix} $P(x,\rho):=\Psi(x,\rho)\tilde\Psi^{-1}(x,\rho)$, where:
$$
\Psi(x,\rho):= \left[
\begin{array}{ll}
\psi_{kk}(x,\rho) & f_{k}(x,\rho)\\
\psi'_{kk}(x,\rho) & f'_k(x,\rho)
\end{array}
\right], \ \rho\in\Omega_+, \
\Psi(x,\rho):= \Psi(x,-\rho), \ \rho\in\Omega_-:=\{\rho: \mbox{Im}\rho<0\}. \eqno(4.1)
$$

\medskip
{\bf Lemma 4.1.} For each fixed $x>0$ $P(x,\rho)$ is bounded for $\rho\to 0$ and $\rho\to\infty$. Moreover, for $\rho\to\infty$ one has $P_1(x,\rho)=I_1+O(\rho^{-1})$, where (and everywhere below) $I$ denotes the identity matrix.

\medskip
{\bf Proof.} By virtue of symmetry it is sufficient to consider $\rho$ from upper half-plane $\overline\Omega_+$. For nonzero $\rho\in\overline\Omega_+$ one has $\det\Psi=\det\tilde\Psi=2i\rho$ and:
$$
2i\rho P_{\xi 1}(x,\rho)=\psi^{(\xi-1)}_{kk}(x,\rho)\tilde f'_k(x,\rho)-f^{(\xi-1)}_k(x,\rho)\tilde\psi'_{kk}(x,\rho),
$$
$$
2i\rho P_{\xi 2}(x,\rho)=f^{(\xi-1)}_k(x,\rho)\tilde\psi_{kk}(x,\rho)-\psi^{(\xi-1)}_{kk}(x,\rho)\tilde f_k(x,\rho).
$$
Lemma 3.2 together with the estimate $f^{(\xi-1)}_k(x,\rho)=O\left(\left|\rho^{\mu_{k1}}\right|\right)$ provide the boundedness of $P(x,\rho)$ for $\rho\to 0$. For $\rho\to\infty$ we use the asymptotics:
$$
f^{(\xi)}_k(x,\rho)=(i\rho)^\xi \exp(i\rho x)[1], \
S^{(\xi)}_{k2}(x,\lambda)=\beta_{k0}\rho^{-\mu_{k2}}\left((-i\rho)^\xi\exp(-i\rho x)[1]+(i\rho)^\xi \gamma_{k0}\exp(i\rho x)[1]\right), \eqno(4.2)
$$
where $[1]:=1+O(\rho^{-1})$ and the constants $\beta_{k0}$, $\gamma_{k0}$ depend only on $\nu_{k0}$ (and therefore are the same for $L$ and $\tilde L$). In particular this yields the estimates:
$$
\hat f^{(\xi)}_k(x,\rho)=O\left((i\rho)^{\xi-1} \exp(i\rho x)\right), \
\hat S^{(\xi)}_{k2}(x,\lambda)=O\left(\rho^{-\mu_{k2}+\xi-1}\exp(-i\rho x)\right). \eqno(4.3)
$$
Moreover, lemma 3.1 provides the estimates:
$$
\hat\Delta(\rho)=O\left(\rho^{\mu_1+2\tau_1-1}\right), \ \Delta^{-1}(\rho)=O\left(\rho^{-\mu_1-2\tau_1}\right). \eqno(4.4)
$$
On the other hand one has the asymptotics (following directly from (3.5)):
$$
\Delta_k(\rho)=\rho^{\mu_1-\mu_{k1}}(d^\infty_k+O(\rho^{-1})), \ d^\infty_k=\sigma^2_k\prod\limits_{j\neq k}\sigma_j b^\infty_{j1}.
$$
Together with lemma 2.1 and the estimates (4.4) this yields:
$$
\gamma_{kk}(\rho)=O\left(\rho^{2\nu_{k1}-2\tau_1}\right)=O(1), \ \hat\gamma_{kk}(\rho)=O(\rho^{-1}\gamma_{kk}(\rho))=O(\rho^{-1}). \eqno(4.5)
$$
The asymptotics (4.2), (4.3), (4.5) provide the estimates:
$$
\psi^{(\xi)}_{kk}(x,\rho)=O\left(\rho^\xi \exp(-i\rho x)\right), \ \hat\psi^{(\xi)}_{kk}(x,\rho)=O\left(\rho^{\xi-1} \exp(-i\rho x)\right). \eqno(4.6)
$$
Together with (4.2), (4.3) the estimates (4.6) yield the required boundedness of $P(x,\rho)$ for $\rho\to\infty$ and the estimate $P_{12}(x,\rho)=O(\rho^{-1})$. In order to obtain the estimate for $P_{11}(x,\rho)-1$ we write:
$$
2i\rho P_{1 1}(x,\rho)=\psi_{kk}(x,\rho)\tilde f'_k(x,\rho)-f_k(x,\rho)\tilde\psi'_{kk}(x,\rho)=
2 i\rho+\hat\psi_{kk}(x,\rho)\tilde f'_k(x,\rho)-\hat f_k(x,\rho)\tilde\psi'_{kk}(x,\rho)
$$
(where we take into account that $\langle\tilde\psi_{kk},\tilde f_k\rangle=2i\rho$) and use again (4.2), (4.3), (4.6). $\hfil\Box$

\medskip
It is clear that $P(x,\rho)$ is meromorphic in $\rho$ outside the real axis and (under condition $G$) for all nonzero real $\rho$ has the limits $P_\pm(x,\rho)$ which are continuous and bounded on $\mathbb{R}\setminus\{0\}$. All the possible poles of $P(x,\rho)$ are necessarily belong to the set $\check{Z}_k:=Z_k\cup\tilde Z_k$, where $Z_k=\{\pm\rho, \rho\in Z_k^+\}$.

\medskip
{\bf Lemma 4.2.} Any $\rho_0\in \check Z_k$ is either a simple pole or a removable singularity for $P(x,\rho)$. Moreover, the following representation holds:
$$
P_{\langle -1\rangle}(x,\rho_0)=\Psi_{\langle 0\rangle}(x,\rho_0)\hat v(\rho_0) (\tilde\Psi^{-1})_{\langle 0\rangle}(x,\rho_0),
$$
where $v(\rho_0):=0$ if $\rho_0\notin Z_k$ and
$$v(\rho_0)=\left[
\begin{array}{cc}
0 & \alpha_k(\rho_0)\\
0 & 0
\end{array}
\right] $$
if $\rho_0\in Z_k$. Here $\alpha_k(\rho_0), \rho_0\in Z_k^+$ are the constants from (3.7) and $\alpha_k(\rho_0):=-\alpha_k(-\rho_0)$ for $\rho_0\in Z_k \cap\Omega_-$.

\medskip
{\bf Proof.}  It follows from (3.7) that $\rho_0\in \check Z_k$ is (at most) a simple pole for $\Psi(x,\rho)$ and the following relation holds:
$$
\Psi_{\langle-1\rangle}(x,\rho_0)=\Psi_{\langle 0\rangle}(x,\rho_0)v(\rho_0). \eqno(4.7)
$$
Moreover, since $\det\Psi=\det\tilde\Psi=\pm 2 i\rho$ for $\pm\rho\in\Omega_+$, the matrix $\tilde\Psi^{-1}(x,\rho)$ has (at most) a simple pole in $\rho_0$ and from $\tilde\Psi\tilde\Psi^{-1}=I$ one can easily obtain:
$$
(\tilde\Psi^{-1})_{\langle-1\rangle}(x,\rho_0)=-\tilde v(\rho_0)(\tilde\Psi^{-1})_{\langle 0\rangle}(x,\rho_0). \eqno(4.8)
$$
Thus, $P(x,\rho)$ has in $\rho_0$ a pole of multiplicity 2 or less. But $P_{\langle -2\rangle}=\Psi_{\langle -1\rangle}(\tilde\Psi^{-1})_{\langle -1\rangle}=-\Psi_{\langle 0\rangle}v\tilde v(\tilde\Psi^{-1})_{\langle 0\rangle}$. We notice that because of special structure of matrices $v(\rho_0)$ and $\tilde v(\rho_0)$ we have  $v(\rho_0)\tilde v(\rho_0)=0$ and therefore $P_{\langle -2\rangle}(x,\rho_0)=0$. Now we calculate $P_{\langle -1\rangle}=\Psi_{\langle -1\rangle}(\tilde\Psi^{-1})_{\langle 0\rangle}+\Psi_{\langle 0\rangle}(\tilde\Psi^{-1})_{\langle -1\rangle}$ and using (4.7) and (4.8) obtain the required representation. $\hfil\Box$

\medskip
Our first result in this section is the following uniqueness theorem.

\medskip
{\bf Theorem 4.1.} $J_k=\tilde J_k$ implies $q_k(x)=\tilde q_k(x)$ for a.e. $x>0$, i.e. specification of the scattering data associated with some ray $\mathcal{R}_k$ uniquely determines the potential $q_k$ on this ray.

\medskip
{\bf Proof.} First we notice that (3.8) yields the following  relation for the matrices $\Psi_\pm(x,\rho)$:
$$
\Psi_+(x,\rho)=\Psi_-(x,\rho)v(\rho), \ \rho\in\mathbb{R}\setminus\{0\} \eqno(4.9)
$$
with "jump matrix"
$$
v(\rho)=\left[
\begin{array}{cc}
r_k(\rho) & 1\\
1-r_k(\rho)r_k(-\rho) & -r_k(-\rho)
\end{array}
\right].  \eqno(4.10)
$$
From (4.10) one can notice that $\tilde r_k=r_k$ implies $\tilde v(\rho)=v(\rho)$ for all $\rho\in\mathbb{R}\setminus\{0\}$ and therefore $P_+(x,\rho)=P_-(x,\rho)$ for each fixed $x>0$ and all nonzero real $\rho$. Thus, if $\tilde J_k=J_k$ the matrix $P(x,\rho)$ is meromorphic in $\mathbb{C}\setminus\{0\}$ with possible poles in points of the set $\check Z_k$. Moreover, $\tilde J_k=J_k$ implies $\check Z_k=\tilde Z_k=Z_k$ and $\tilde v(\rho_0)=v(\rho_0)$ for any $\rho_0\in Z_k$. Then, by virtue of lemma 4.2 we have $P_{\langle -1\rangle}(x,\rho_0)=0$, i.e. $P(x,\rho)$ has actually  removable singularities in all $\rho_0\in Z_k$. Thus, $P(x,\rho)$ is holomorphic in $\mathbb{C}\setminus\{0\}$ and by virtue of lemma 4.1 bounded for $\rho\to 0$ and $\rho\to \infty$. Moreover, by virtue of lemma 4.1 $P_1(x,\rho)-I_1$ vanishes at infinity. This means that $P_1(x,\rho)-I_1$ is actually equal to 0 identically. Thus, we have $P_1(x,\rho)\equiv I_1$ that yields $f_k(x,\rho)\equiv\tilde f_k(x,\rho)$ and subsequently $q_k(x)=\tilde q_k(x)$ for a.e. $x>0$. $\hfil\Box$

\medskip
Our next goal is a constructive procedure for solving the problem $IP(k)$. From this point we assume that the numbers $\nu_{j0}$, $j=\overline{1,p}$ and coefficients of the forms $U_{j1}$, $U_{j2}$ in matching conditions (3.2) are known. Also we suppose that $\tilde L$ is chosen a priori "model" (or "reference") problem with known potentials $q_j, j=\overline{1,p}$. Below we reduce the problem $IP(k)$ to a certain linear equation ("main equation") in some Banach space and prove the unique solvability of this equation.

We consider again the spectral mappings matrix $P(x,\rho)$. By virtue of lemma 4.1 one has $\int_{|\mu|=R}(\mu-\rho)^{-1}d\mu (P_1(x,\mu)-I_1) \to 0$ as $R\to\infty$  that yields the following basic relation:
$$
P_1(x,\rho)-I_1=\sum\limits_{\mu\in \check Z_k}(\rho-\mu)^{-1}P_{1,\langle-1\rangle}(x,\mu)+ \frac{1}{2\pi i}\int\limits_{-\infty}^\infty  \frac{d\mu}{\mu-\rho}\left(P_{1,+}(x,\mu)-P_{1,-}(x,\mu)\right) \eqno(4.11)
$$
with arbitrary $\rho\in\mathbb{C}\setminus(\mathbb{R}\cup\check Z_k)$.

Let us introduce the matrix $V(x,\rho):=P^{-1}_-(x,\rho)P_+(x,\rho)=$  $\tilde\Psi_-(x,\rho)v(\rho)\tilde\Psi^{-1}_+(x,\rho)$, $\rho\in\mathbb{R}\setminus\{0\}$. We notice that $V$ is uniquely determined by $\tilde L$ and scattering data $J_k$. We also note that in view of lemma 4.1 and $\det P=1$, $V(x,\rho)$ is continuous and bounded with respect to $\rho\in\mathbb{R}\setminus\{0\}$ and any fixed $x>0$.

We define
$$
d_1(x,\rho,\mu_0)=\left.\left[(\rho-\mu)^{-1} \Psi^{-1}(x,\mu)\right]_{\langle 0\rangle}\right|_{\mu=\mu_0}, \ \rho\in\mathbb{C}\setminus\check Z_k, \mu_0\in \check Z_k $$$$
d_2(x,\rho_0,\mu)=\left.\left[(\rho-\mu)^{-1} \Psi(x,\rho)\right]_{\langle 0\rangle}\right|_{\rho=\rho_0}, \ \mu\in\mathbb{C}\setminus\check Z_k, \rho_0\in \check Z_k $$$$
D(x,\rho,\mu):=(\rho-\mu)^{-1}( \Psi^{-1})(x,\mu)\Psi(x,\rho), \ \rho\neq\mu, \rho,\mu\in\mathbb{C}\setminus\check Z_k, $$$$
D(x,\rho,\mu_0):=\left.\left[D(x,\rho,\mu)\right]_{\langle 0\rangle}\right|_{\mu=\mu_0}, \ \rho\in\mathbb{C}\setminus\check Z_k, \mu_0\in \check Z_k $$$$
D(x,\rho_0,\mu):=\left.\left[D(x,\rho_0,\mu)\right]_{\langle 0\rangle}\right|_{\rho=\rho_0}, \ \mu\in\mathbb{C}\setminus\check Z_k, \rho_0\in \check Z_k
$$
$$
D(x,\rho_0,\mu_0):=\left.\left[D(x,\rho,\mu_0)\right]_{\langle 0\rangle}\right|_{\rho=\rho_0}=\left.\left[D(x,\rho_0,\mu)\right]_{\langle 0\rangle}\right|_{\mu=\mu_0},\ \rho_0\in \check Z_k,\ \mu_0\in \check Z_k
$$
and $\tilde d_{j}(x,\rho,\mu), \ j=1,2$, $\tilde D(x,\rho,\mu)$ by similar formulae, where $\Psi(x,\rho)$ is replaced with $\tilde \Psi(x,\rho)$. We also define:
$$
\tilde A(x,\rho,\mu):= \tilde D(x,\rho,\mu)\hat v(\rho), \  A(x,\rho,\mu):= - D(x,\rho,\mu)\hat v(\rho), \ \rho\in \check Z_k,\ \mu\in \check Z_k.
$$

\medskip
{\bf Lemma 4.3.} For $\mu\in\check Z_k$, $\rho\notin\check Z_k$ one has:
$$
(\rho-\mu)^{-1}P_{\langle -1\rangle}(x,\mu)=\Psi_{\langle 0\rangle}(x,\mu)\hat v(\mu) \tilde d_1(x,\rho,\mu),
$$
$$
(\rho-\mu)^{-1}P_{\langle -1\rangle}(x,\mu)\tilde\Psi(x,\rho)=\Psi_{\langle 0\rangle}(x,\mu)\hat v(\mu) \tilde D(x,\rho,\mu).
$$
For $\xi\in\check Z_k$, $\rho,\mu\notin\check Z_k$
$$
(\rho-\xi)^{-1}(\xi-\mu)^{-1}\Psi^{-1}(x,\mu)P_{\langle -1\rangle}(x,\xi)\tilde\Psi(x,\rho)=D(x,\xi,\mu)\hat v(\xi) \tilde D(x,\rho,\xi).
$$
\medskip
{\bf Proof.} The relations can be easily verified by direct calculation. Consider, for instance,
$$
\tilde d_1(x,\rho,\mu_0)=\left.\left[(\rho-\mu)^{-1} \tilde\Psi^{-1}(x,\mu)\right]_{\langle 0\rangle}\right|_{\mu=\mu_0}=
$$
$$
(\rho-\mu_0)^{-1} (\tilde\Psi^{-1})_{\langle 0\rangle}(x,\mu_0)+
\left.\left[(\rho-\mu)^{-1}\right]_{\langle 1\rangle}\right|_{\mu=\mu_0}\cdot(\tilde\Psi^{-1})_{\langle -1\rangle}(x,\mu_0).
$$
Multiplying this by $\Psi_{\langle 0\rangle}(x,\mu_0)\hat v(\mu_0)$, then using lemma 4.2, relation (4.8) and taking into account that $\hat v(\mu_0)\tilde v(\mu_0)=0$ we obtain the first of the required relations. The others can be obtained in a similar way.$\hfil\Box$

\medskip
The assertion below follows directly from definition of $d_j(x,\rho,\mu)$.

\medskip
{\bf Lemma 4.4.} For any fixed  $\mu\in\check Z_k$ $d_1(x,\rho,\mu)$ is holomorphic in $\rho\in\mathbb{C}\setminus\check Z_k$. In particular for real $\rho$ one has $d_1(x,\rho+i0,\mu)=d_1(x,\rho-i0,\mu)=d_1(x,\rho,\mu)$, moreover, $d_1(x,\cdot,\mu)\in L_2\left(\mathbb{R}, \mathbb{C}^2\right)$. Analogously,  for any fixed  $\rho\in\check Z_k$ $d_2(x,\rho,\mu)$ is holomorphic in $\mu\in\mathbb{C}\setminus\check Z_k$, $d_2(x,\rho,\cdot)\in L_2\left(\mathbb{R}, \mathbb{C}^2\right)$. Here and below elements of $\mathbb{C}^2$ are considered as row-vectors.

\medskip
In our further considerations we will also use the Cauchy operators:
$$
(Cf)(\rho):=\frac{1}{2\pi i}\int\limits_{-\infty}^\infty \frac{d\mu}{\mu-\rho}f(\mu), \ \rho\in\mathbb{C}\setminus\mathbb{R},
$$
$$
(C_\pm f)(\rho):=(Cf)_\pm(\rho), \ \rho\in\mathbb{R}
$$
and the Plemelj-Sokhotskii formula $C_\pm=\mathcal{P}\pm E/2$. Here (and everywhere below) $E$ denotes the identical operator and
$$(\mathcal{P}f)(\rho)=\frac{1}{2\pi i}\int\limits_{-\infty}^\infty \frac{d\mu}{\mu-\rho}f(\mu), \ \rho\in\mathbb{R}$$
where a principal value of integral is assumed. We recall that $C_\pm$ and $\mathcal{P}$ are bounded operators in $L_2(\mathbb{B}, \mathcal{M})$ (for any finite-dimensional space $\mathcal{M})$ and the following relation holds for two arbitrary matrices $F_1, F_2$ of suitable orders with quadratically integrable elements:
$$
\int\limits_{-\infty}^\infty d\mu (C_\pm F_1)(\mu) F_2(\mu)=-\int\limits_{-\infty}^\infty d\mu F_1(\mu) (C_\mp F_2)(\mu).
$$

\medskip
Define $\Phi(x,\rho):=P_{1,+}(x,\rho)-P_{1,-}(x,\rho)$ for real $\rho$ and $\Phi(x,\rho):=\Psi_{1,\langle 0\rangle}(x,\rho)\hat v(\rho)$ for $\rho\in\check Z_k$. One can easily notice that by virtue of lemma 4.1 $\left.\Phi(x,\cdot)\right|_\mathbb{R}\in L_2\left(\mathbb{R},\mathbb{C}^2\right)$. Now we return to the relation (4.11) and using the notations introduced above and lemma 4.3  rewrite it as follows:
$$
P_1(x,\rho)-I_1=\sum\limits_{\mu\in \check Z_k}\Phi(x,\mu)\tilde d_1(x,\rho,\mu)+(C\Phi(x,\cdot))(\rho), \ \rho\in\mathbb{C}\setminus(\mathbb{R}\cup\check Z_k). \eqno(4.12)
$$
Taking in (4.12) the limits as $\pm\mbox{Im}\rho\to 0$ and substituting them into the jump relation $P_{1,+}-P_{1,-}V=0$ we arrive at:
$$
\Phi(x,\rho)\left(\frac{1}{2}I+\frac{1}{2}V(x,\rho)\right)+
\frac{1}{2\pi i}\int\limits_{-\infty}^\infty \frac{d\mu}{\mu-\rho}\Phi(x,\mu)(I-V(x,\rho))+$$$$
\sum\limits_{\mu\in \check Z_k}\Phi(x,\mu)\tilde d_1(x,\rho,\mu)(I-V(x,\rho))+I_1-V_1(x,\rho)=0, \ \rho\in\mathbb{R}\setminus\{0\}, \eqno(4.13)
$$
where a principal value of integral is assumed. In the particular case when $\check Z_k=\emptyset$ the relation (4.13)  can be considered for each fixed $x>0$ as a linear equation with respect to $\Phi(x,\cdot)$. In general case (4.13) should be completed with some relations at the points $\rho\in\check Z_k$. We proceed as follows. Multiplying (4.12) by $\tilde\Psi(x,\rho)$ and using lemma 4.3  we obtain:
$$
\Psi_{1}(x,\rho)-\tilde\Psi_{1}(x,\rho)=\sum\limits_{\mu\in \check Z_k} \Psi_{1,\langle 0\rangle}(x,\mu)\hat v(\mu)\tilde D(x,\rho,\mu)+\left(C\Phi(x)\right)(\rho)\tilde\Psi(x,\rho), \ \rho\in\mathbb{C}\setminus(\mathbb{R}\cup\check Z_k).
$$
Now we take an arbitrary $\rho_0\in \check Z_k$ and multiply the relation above by $\hat v(\rho_0)$. Thus we get:
$$
\Psi_{1}(x,\rho)\hat v(\rho_0)-\tilde\Psi_{1}(x,\rho)\hat v(\rho_0)=$$$$\sum\limits_{\mu\in \check Z_k} \Psi_{1,\langle 0\rangle}(x,\mu)\hat v(\mu)\tilde D(x,\rho,\mu)\hat v(\rho_0)+\frac{1}{2\pi i}\int\limits_{-\infty}^\infty \frac{d\mu}{\mu-\rho}\Phi(x,\mu) \tilde\Psi(x,\rho)\hat v(\rho_0). \eqno(4.14)
$$
Consider again the relation (4.14). Taking the coefficient $\left.[\dots]_{\langle 0 \rangle}\right|_{\rho=\rho_0}$ in Laurent series of its both sides and using lemma 4.3 we arrive at:
$$
\Phi(x,\rho_0)+\frac{1}{2\pi i}\int\limits_{-\infty}^\infty d\mu\Phi(x,\mu) \tilde d_2(x,\rho_0,\mu)\hat v(\rho_0)-\sum\limits_{\mu\in \check Z_k} \Phi(x,\mu)\tilde A(x,\rho_0,\mu)=\tilde\Psi_{1,\langle 0 \rangle}(x,\rho_0)\hat v(\rho_0), \rho_0\in\check Z_k, \eqno(4.15)
$$
where $\rho_0\in\check Z_k$ is arbitrary. Together with (4.13) (4.15) forms the complete linear system with respect to $\Phi(x,\rho)$, $\rho\in\mathbb{R}\cup\check Z_k$. We can summarize our considerations as follows.

\medskip
{\bf Theorem 4.2.} Let $\tilde d_j(x,\rho,\mu)$ and $\tilde A(x,\rho,\mu)$ constructed as described above by given scattering data $J_k$ and given model problem $\tilde L$. Then
$$
\Phi(x,\rho)=\left\{\begin{array}{cc}
P_{1,+}(x,\rho)-P_{1,-}(x,\rho), \ \rho\in\mathbb{R}, \\
\Psi_{1,\langle 0\rangle}(x,\rho)\hat v(\rho), \rho\in\check Z_k
\end{array}
\right.
$$
for each fixed $x>0$ solves the linear system (4.13), (4.15).

\medskip
In order to complete our procedure of solving the problem $IP(k)$ we are to show the unique solvability of the specified linear system. First we rewrite the system as a linear equation in Banach space $\mathcal{H}:=\mathcal{H}_r\oplus\mathcal{H}_d$, where $\mathcal{H}_r=L_2\left(\mathbb{R}, \mathbb{C}^2\right)$, $\mathcal{H}_d=\left(\mathbb{C}^2\right)^{\check Z_k}$. As we already mentioned above for each fixed $x>0$ $\Phi(x,\cdot)$ can be considered as an element of $\mathcal{H}$. Thus, the system (4.13), (4.15) can be written in the form $\mathcal{A}\Phi=G$, where
$$
G(\rho):=\left\{\begin{array}{cc}
V_1(x,\rho)-I_1, \ \rho\in\mathbb{R}, \\
\tilde\Psi_{1,\langle 0\rangle}(x,\rho)\hat v(\rho), \rho\in\check Z_k,
\end{array}
\right.
$$
while $\mathcal{A}$ is a linear bounded operator in $\mathcal{H}$ defined by the following operator matrix:
$$
\mathcal{A} =\left[
\begin{array}{ll}
\mathcal{A}_{rr}  & \mathcal{A}_{rd}  \\
\mathcal{A}_{dr}  & \mathcal{A}_{dd}
\end{array}
\right],
$$
$$
\mathcal{A}_{rr}  \varphi=C_+\vp - (C_-\vp)V,\
(\mathcal{A}_{rd}  \varphi)(\rho)=\sum\limits_{\mu\in \check Z_k}\varphi(\mu)\tilde d_1(\rho,\mu)(I-V(\rho)), \ \rho\in\mathbb{R},
$$
$$
(\mathcal{A}_{dr}  \varphi)(\rho)= \frac{1}{2\pi i}\int\limits_{-\infty}^\infty \vp(\xi)\tilde d_2(\rho,\xi)\hat v(\rho)d\xi, \  \rho\in\check Z_k ,\
(\mathcal{A}_{dd}  \varphi)(\rho)=\vp(\rho)-\sum\limits_{\mu\in \check Z_k}\varphi(\mu)\tilde A(x,\rho,\mu), \ \rho\in\check Z_k .
$$
Here and below we assume that parameter $x>0$ is (arbitrary) fixed and for the sake of brevity we omit it in all the argument's lists. Boundedness of the operator $\mathcal{A}$ follows from lemma 4.4 and boundedness of operators $C_\pm$, function $V$.
 Now we can provide the following improvement of Theorem 4.2.

\medskip
{\bf Theorem 4.3.} For each fixed $x>0$ $\Phi(x,\cdot)$ is a unique solution in $\mathcal{H}$ of the equation $\mathcal{A}\phi=G$. Moreover, operator $\mathcal{A}$ has a bounded inverse operator $\mathcal{A}^{-1}=\mathcal{B}$, where:
$$
\mathcal{B}=\left[
\begin{array}{ll}
 \mathcal{B}_{rr}  &  \mathcal{B}_{rd}  \\
 \mathcal{B}_{dr}  &  \mathcal{B}_{dd}
\end{array}
\right],$$
$$
\mathcal{B}_{rr}  f =C_+(f\tilde P_+)P_+-C_-(f\tilde P_+)P_-, \
( \mathcal{B}_{rd}  f)(\rho)=\sum\limits_{\mu\in \check Z_k}f(\mu) d_1(\rho,\mu)(P_+(\rho)-P_-(\rho)), \ \rho\in\mathbb{R},
$$
$$
(\mathcal{B}_{dr}  f)(\rho)= -\frac{1}{2\pi i}\int\limits_{-\infty}^\infty  f(\xi)\tilde P_+(\xi) d_2(\rho,\xi)\hat v(\rho)d\xi, \  \rho\in\check Z_k, \
(\mathcal{B}_{dd}  f)(\rho)=f(\rho)-\sum\limits_{\mu\in \check Z_k}f(\mu) A(x,\rho,\mu), \ \rho\in\check Z_k.
$$
Here and below $\tilde P:=P^{-1}$ while $P$  denotes as above the spectral mappings matrix.

\medskip
{\bf Lemma 4.5.} The following relation holds for nonreal $\rho\neq\mu$, $\rho,\mu\in\mathbb{C} \setminus\check Z_k$:
$$
\frac{1}{2\pi i}\int\limits_{-\infty}^\infty \frac{d\xi}{(\rho-\xi)(\xi-\mu)}(P_+(\xi)-P_-(\xi))
=\frac{1}{\rho-\mu}P(\mu)-\frac{1}{\rho-\mu}P(\rho)+\sum\limits_{\xi\in \check Z_k} d_2(\xi,\mu)\hat v(\xi) \tilde d_1(\rho,\xi).
$$
For $\rho\in\mathbb{C}\setminus(\mathbb{R}\cup \check Z_k),\mu\in \check Z_k$ one has:
$$
\frac{1}{2\pi i}\int\limits_{-\infty}^\infty \frac{d\xi}{\rho-\xi}d_1(\xi,\mu)(P_+(\xi)-P_-(\xi))=\tilde d_1(\rho,\mu)-d_1(\rho,\mu)P(\rho)
-\sum\limits_{\xi\in \check Z_k}A(\xi,\mu)\tilde d_1(\rho,\xi).
$$
For $\mu\in\mathbb{C}\setminus(\mathbb{R}\cup \check Z_k),\rho\in \check Z_k$ one has:
$$
\frac{1}{2\pi i}\int\limits_{-\infty}^\infty \frac{d\xi}{\xi-\mu}(P_+(\xi)-P_-(\xi))\tilde d_2(\rho,\xi)\hat v(\rho)=P(\mu)\tilde d_2(\rho,\mu)\hat v(\rho)-
 d_2(\rho,\mu)\hat v(\rho)+\sum\limits_{\xi\in \check Z_k} d_2(\xi,\mu)\hat v(\xi)\tilde A(\rho,\xi).
$$
For $\rho,\mu\in \check Z_k$ one has:
$$
\frac{1}{2\pi i}\int\limits_{-\infty}^\infty d_1(\xi,\mu)(P_+(\xi)-P_-(\xi))\tilde d_2(\rho,\xi)\hat v(\rho)d\xi=
\tilde A(\rho,\mu)+A(\rho,\mu)-\sum\limits_{\xi\in \check Z_k} A(\xi,\mu)\tilde A(\rho,\xi).$$
Symmetrically, one can obtain the following relations:
$$
\frac{1}{2\pi i}\int\limits_{-\infty}^\infty \frac{d\xi}{(\rho-\xi)(\xi-\mu)}(\tilde P_+(\xi)-\tilde P_-(\xi))
=\frac{1}{\rho-\mu}\tilde P(\mu)-\frac{1}{\rho-\mu}\tilde P(\rho)-\sum\limits_{\xi\in \check Z_k} \tilde d_2(\xi,\mu)\hat v(\xi)  d_1(\rho,\xi)
$$
for $\rho,\mu\in\mathbb{C}\setminus(\mathbb{R}\cup \check Z_k)$,
$$
\frac{1}{2\pi i}\int\limits_{-\infty}^\infty \frac{d\xi}{\rho-\xi}\tilde d_1(\xi,\mu)(\tilde P_+(\xi)-\tilde P_-(\xi))= d_1(\rho,\mu)-\tilde d_1(\rho,\mu)\tilde P(\rho)
-\sum\limits_{\xi\in \check Z_k}\tilde A(\xi,\mu) d_1(\rho,\xi)$$
for $\rho\in\mathbb{C}\setminus(\mathbb{R}\cup \check Z_k),\mu\in \check Z_k$,
$$
\frac{1}{2\pi i}\int\limits_{-\infty}^\infty \frac{d\xi}{\xi-\mu}(\tilde P_+(\xi)-\tilde P_-(\xi)) d_2(\rho,\xi)\hat v(\rho)=\tilde P(\mu) d_2(\rho,\mu)\hat v(\rho)-
 \tilde d_2(\rho,\mu)\hat v(\rho)+\sum\limits_{\xi\in \check Z_k}\tilde d_2(\xi,\mu)\hat v(\xi) A(\rho,\xi)
$$
for $\mu\in\mathbb{C}\setminus(\mathbb{R}\cup \check Z_k),\rho\in \check Z_k$,
$$
-\frac{1}{2\pi i}\int\limits_{-\infty}^\infty \tilde d_1(\xi,\mu)(\tilde P_+(\xi)-\tilde P_-(\xi)) d_2(\rho,\xi)\hat v(\rho)d\xi=
\tilde A(\rho,\mu)+A(\rho,\mu)-\sum\limits_{\xi\in \check Z_k}\tilde A(\xi,\mu) A(\rho,\xi)$$
for $\rho,\mu\in \check Z_k$.

\medskip
{\bf Proof.} All the relations can be obtained in a similar way based on the relation:
$$
\lim\limits_{R\to\infty}\int\limits_{|\xi|=R}\frac{d\xi}{(\rho-\xi)(\xi-\mu)}P(\xi)=0, \ \rho,\mu\in\mathbb{C}\setminus(\mathbb{R}\cup\check Z_k), \eqno(4.16)
$$
that follows from lemma 4.1.

Let us show how to get, for instance, the forth relation. Multiplying (4.16) by $\Psi^{-1}(\mu)$ from the left, by $\tilde\Psi(\rho)$ from the right and applying the residue theorem we obtain:
$$
\frac{1}{2\pi i}\int\limits_{-\infty}^\infty\frac{d\xi}{(\rho-\xi)(\xi-\mu)}\Psi^{-1}(\mu)(P_+(\xi)-P_-(\xi))\tilde\Psi(\rho)=$$$$
\sum\limits_{\xi\in\check Z_k}\frac{1}{(\rho-\xi)(\xi-\mu)}\Psi^{-1}(\mu)P_{\langle -1\rangle}(\xi)\tilde\Psi(\rho)+
\frac{1}{\rho-\mu}\tilde\Psi^{-1}(\mu)\tilde\Psi(\rho)-\frac{1}{\rho-\mu}\Psi^{-1}(\mu)\Psi(\rho), \ \rho,\mu\in\mathbb{C}\setminus(\mathbb{R}\cup\check Z_k).
$$
Using lemma 4.3 we rewrite this as follows:
$$
\frac{1}{2\pi i}\int\limits_{-\infty}^\infty\frac{d\xi}{(\rho-\xi)(\xi-\mu)}\Psi^{-1}(\mu)(P_+(\xi)-P_-(\xi))\tilde\Psi(\rho)=$$$$
\sum\limits_{\xi\in\check Z_k}D(\xi,\mu)\hat v(\xi)\tilde D(\rho,\xi)+
\tilde D(\rho,\mu)-D(\rho,\mu), \ \rho,\mu\in\mathbb{C}\setminus(\mathbb{R}\cup\check Z_k).
$$
Taking (for arbitrary fixed $\rho\in\mathbb{C}\setminus(\mathbb{R}\cup\check Z_k)$, $\mu_0\in\check Z_k$) the coefficients $\left.[...]_{\langle 0\rangle}\right|_{\mu=\mu_0}$ in Laurent series of both sides we arrive at:
$$
\frac{1}{2\pi i}\int\limits_{-\infty}^\infty\frac{d\xi}{\rho-\xi}d_1(\xi,\mu_0)(P_+(\xi)-P_-(\xi))\tilde\Psi(\rho)=$$$$-\sum\limits_{\xi\in\check Z_k}A(\xi,\mu_0)\tilde D(\rho,\xi)+\tilde D(\rho,\mu_0)-D(\rho,\mu_0), \ \mu_0\in\check Z_k, \rho\in\mathbb{C}\setminus(\mathbb{R}\cup\check Z_k).
$$
Finally, taking (for arbitrary $\rho_0,\mu_0\in\check Z_k$) the coefficients $\left.[...]_{\langle 0\rangle}\right|_{\rho=\rho_0}$ in Laurent series of both sides and multiplying them by $\hat v(\rho_0)$ we obtain the required relation. $\hfil\Box$

\medskip
{\bf Proof of Theorem 4.3.}
Boundedness of the operator $\mathcal{B}$ follows from lemma 4.4 and boundedness of operators $C_\pm$, functions $V, P_\pm, \tilde P_\pm$. The proof consists in the direct calculations of elements of the operator matrices $\mathcal{AB}$ and $\mathcal{BA}$.

{\bf Position (r,r).} Let $\vp=\mathcal{B}_{rr} f$. We rewrite it as $\vp(\rho)=f(\rho)+\left[C_-(f\tilde P_+)\right](\rho)(P_+(\rho)-P_-(\rho))$. Then $(C\vp)(\rho)$ for arbitrary fixed nonreal $\rho\in\mathbb{C}\setminus\check Z_k$ can be written in the following form:
$$
(C\vp)(\rho)=(Cf)(\rho)-\frac{1}{2\pi i}\int\limits_{-\infty}^\infty  \left[C_-(f\tilde P_+)\right](\mu) F_\rho(\mu)d\mu,
$$
where $F_\rho(\mu):=(\rho-\mu)^{-1}(P_+(\mu)-P_-(\mu))$. Using the relation
$$
-\int\limits_{-\infty}^\infty  (C_-F_1)(\mu) F_2(\mu)d\mu=\int\limits_{-\infty}^\infty  F_1(\mu) (C_+F_2)(\mu)d\mu
$$
we rewrite it as follows:
$$
(C\vp)(\rho)=(Cf)(\rho)+\frac{1}{2\pi i}\int\limits_{-\infty}^\infty  (f\tilde P_+)(\mu) (C_+F_\rho)(\mu)d\mu.
$$
Now let us consider $(CF_\rho)(\mu)$. For nonreal $\mu\in \mathbb{C}\setminus\check Z_k$ using lemma 4.5 we calculate:
$$
(CF_\rho)(\mu)=\frac{1}{2\pi i}\int\limits_{-\infty}^\infty \frac{d\xi}{(\rho-\xi)(\xi-\mu)}(P_+(\xi)-P_-(\xi))=
\frac{1}{\rho-\mu}P(\mu)-\frac{1}{\rho-\mu}P(\rho)+\sum\limits_{\xi\in \check Z_k} d_2(\xi,\mu)\hat v(\xi) \tilde d_1(\rho,\xi).
$$
Taking the limit as $\mbox{Im}\mu\to +0$ (while $\rho\in \mathbb{C}\setminus(\mathbb{R}\cup\check Z_k)$ remains fixed) we obtain:
$$
(C_+F_\rho)(\mu)=\frac{1}{\rho-\mu}P_+(\mu)-\frac{1}{\rho-\mu}P(\rho)+\sum\limits_{\xi\in \check Z_k} d_2(\xi,\mu)\hat v(\xi) \tilde d_1(\rho,\xi).
$$
Thus, we can write:
$$
(C\vp)(\rho)=(Cf)(\rho)+\frac{1}{2\pi i}\int\limits_{-\infty}^\infty (f\tilde P_+)(\mu)
\left\{\frac{1}{\rho-\mu}P_+(\mu)-\frac{1}{\rho-\mu}P(\rho)+\sum\limits_{\xi\in \check Z_k} d_2(\xi,\mu)\hat v(\xi) \tilde d_1(\rho,\xi)\right\}d\mu=
$$
$$
\frac{1}{2\pi i}\int\limits_{-\infty}^\infty (f\tilde P_+)(\mu)
\left\{\sum\limits_{\xi\in \check Z_k} d_2(\xi,\mu)\hat v(\xi) \tilde d_1(\rho,\xi)\right\}d\mu+\left[C(f\tilde P_+)\right](\rho)P(\rho)
$$
that yields the following relation:
$$
(\mathcal{A}_{rr} \mathcal{B}_{rr} f)(\rho)=(\mathcal{A}_{rr} \vp)(\rho)=\left[C_+(f\tilde P_+)\right](\rho)P_+(\rho)-\left[C_-(f\tilde P_+)\right](\rho)P_-(\rho)V(\rho)+
$$
$$
\frac{1}{2\pi i}\int\limits_{-\infty}^\infty (f\tilde P_+)(\mu)
\left\{\sum\limits_{\xi\in \check Z_k} d_2(\xi,\mu)\hat v(\xi) \tilde d_1(\rho,\xi)\right\}(I-V(\rho))d\mu.
$$
On the other hand we have:
$$
(\mathcal{A}_{rd} \mathcal{B}_{dr} f)(\rho)=\sum\limits_{\xi\in\check Z_k}\left\{-\frac{1}{2\pi i}\int\limits_{-\infty}^\infty  f(\mu)\tilde P_+(\mu) d_2(\xi,\mu)\hat v(\xi)d\mu\right\}\tilde d_1(\rho,\xi)(I-V(\rho))
$$
and thus arrive at:
$$
(\mathcal{A}_{rr} \mathcal{B}_{rr} f)(\rho)+(\mathcal{A}_{rd} \mathcal{B}_{dr} f)(\rho)=\left[C_+(f\tilde P_+)\right](\rho)P_+(\rho)-\left[C_-(f\tilde P_+)\right](\rho)P_-(\rho)V(\rho).
$$
Taking into account that $P_+=P_-V$ and $\tilde P=P^{-1}$ we obtain finally $(\mathcal{A}_{rr} \mathcal{B}_{rr} f)(\rho)+(\mathcal{A}_{rd} \mathcal{B}_{dr} f)(\rho)=[(C_+-C_-)(f\tilde P_+)](\rho)P_+(\rho)=f(\rho)$.

Now let $f=\mathcal{A}_{rr} \vp$. Then $(f\tilde P_+)(\rho)=(\vp\tilde P_+)(\rho)+(C_-\vp)(\rho)(\tilde P_+(\rho)-\tilde P_-(\rho))$. Consider the function $[C(f\tilde P_+)](\rho)$. Proceeding as above we obtain:
$$
[C(f\tilde P_+)](\rho)=[C(\vp\tilde P_+)](\rho)+\frac{1}{2\pi i}\int\limits_{-\infty}^\infty  \vp(\mu)[C_+\tilde F_\rho](\mu)d\mu,
$$
where $\tilde F_\rho(\mu):=(\rho-\mu)^{-1}(\tilde P_+(\mu)-\tilde P_-(\mu))$. Using lemma 4.5 we calculate:
$$
[C_+\tilde F_\rho](\mu)=\frac{1}{\rho-\mu}\tilde P_+(\mu)-\frac{1}{\rho-\mu}\tilde P(\rho)-
\sum\limits_{\xi\in\check Z_k}\tilde d_2(\xi,\mu)\hat v(\xi) d_1(\rho,\xi),
$$
$$
[C(f\tilde P_+)](\rho)P(\rho)=(C\vp)(\rho)-\frac{1}{2\pi i}\int\limits_{-\infty}^\infty \vp(\mu)\left\{\sum\limits_{\xi\in\check Z_k}\tilde d_2(\xi,\mu)\hat v(\xi) d_1(\rho,\xi)\right\} P(\rho)d\mu
$$
that yields finally:
$$
(\mathcal{B}_{rr} \mathcal{A}_{rr} \vp)(\rho)=(\mathcal{B}_{rr} f)(\rho)=\vp(\rho)-\frac{1}{2\pi i}\int\limits_{-\infty}^\infty \vp(\mu)\left\{\sum\limits_{\xi\in\check Z_k}\tilde d_2(\xi,\mu)\hat v(\xi) d_1(\rho,\xi)\right\} (P_+(\rho)-P_-(\rho))d\mu.
$$
On the other hand we have:
$$
(\mathcal{B}_{rd} \mathcal{A}_{dr} \vp)(\rho)=\sum\limits_{\xi\in\check Z_k}
\left\{\frac{1}{2\pi i}\int\limits_{-\infty}^\infty \vp(\mu)\tilde d_2(\xi,\mu)\hat v(\xi) d\mu\right\}d_1(\rho,\xi)(P_+(\rho)-P_-(\rho))
$$
and thus $(\mathcal{B}_{rr} \mathcal{A}_{rr} \vp)(\rho)+(\mathcal{B}_{rd} \mathcal{A}_{dr} \vp)(\rho)=\vp(\rho)$.

{\bf Position (r,d).} Calculate $\mathcal{A}_{rr} \mathcal{B}_{rd} +\mathcal{A}_{rd} \mathcal{B}_{dd} $. Let $\vp=\mathcal{B}_{rd} f$. Then
$$
\mathcal{A}_{rr} \vp(\rho)=\sum\limits_{\mu\in \check Z_k} f(\mu)\left[C_+R_\mu-(C_-R_\mu)V\right](\rho),
$$
where $R_\mu(\rho):=d_1(\rho,\mu)(P_+(\rho)-P_-(\rho))$. Consider for nonreal $\rho\in \mathbb{C}\setminus\check Z_k$ the function:
$$
CR_\mu(\rho)=-\frac{1}{2\pi i}\int\limits_{-\infty}^\infty \frac{d\xi}{\rho-\xi}d_1(\xi,\mu)(P_+(\xi)-P_-(\xi)).
$$
Using lemma 4.5 we rewrite this as follows:
$$
CR_\mu(\rho)=-\tilde d_1(\rho,\mu)+d_1(\rho,\mu)P(\rho)
+\sum\limits_{\xi\in \check Z_k}A(\xi,\mu)\tilde d_1(\rho,\xi).
$$
Taking the limits $\pm\mbox{Im}\rho\to 0$ we obtain:
$$
\left[C_+R_\mu-(C_-R_\mu)V\right](\rho)=-\tilde d_1(\rho,\mu)(I-V(\rho))+d_1(\rho,\mu)(P_+(\rho)-P_-(\rho)V(\rho))$$$$
+\sum\limits_{\xi\in \check Z_k}A(\xi,\mu)\tilde d_1(\rho,\xi)(I-V(\rho)).
$$
Thus, taking into account that $P_+=P_-V$ we arrive at:
$$
(\mathcal{A}_{rr} \mathcal{B}_{rd} )(\rho)=\left\{\sum\limits_{\mu\in\check Z_k}f(\mu)F(\rho,\mu)\right\}(I-V(\rho)),
$$
where
$$
F(\rho,\mu)=\sum\limits_{\xi\in \check Z_k}A(\xi,\mu)\tilde d_1(\rho,\xi)-\tilde d_1(\rho,\mu).
$$
On the other hand direct calculation yields:
$$
(\mathcal{A}_{rd} \mathcal{B}_{dd} f)(\rho)=\left\{\sum\limits_{\xi\in \check Z_k}\left(f(\xi)-\sum\limits_{\mu\in\check Z_k}f(\mu) A(\xi,\mu)\tilde d_1(\rho,\xi)\right)\right\}(I-V(\rho)),
$$
that can be rewritten as
$$
(\mathcal{A}_{rd} \mathcal{B}_{dd} )(\rho)=-\left\{\sum\limits_{\mu\in\check Z_k}f(\mu)F(\rho,\mu)\right\}(I-V(\rho)).
$$
Thus, we have $\mathcal{A}_{rr} \mathcal{B}_{rd} +\mathcal{A}_{rd} \mathcal{B}_{dd} =0$. Symmetrical calculations show that  $\mathcal{B}_{rr} \mathcal{A}_{rd} +\mathcal{B}_{rd} \mathcal{A}_{dd} =0$.

\medskip
{\bf Position (d,r)}.  Calculate $\mathcal{A}_{dr} \mathcal{B}_{rr} +\mathcal{A}_{dd} \mathcal{B}_{dr} $. Let $\vp=\mathcal{B}_{dr} f$. Then $\mathcal{A}_{dd} \vp(\rho)$ can be written in the following form:
$$
\mathcal{A}_{dd} \vp(\rho)=\frac{1}{2\pi i}\int\limits_{-\infty}^\infty f(\mu)\tilde P_+(\mu)
\left[-d_2(\rho,\mu)\hat v(\rho)+\sum\limits_{\xi\in\check Z_k}d_2(\xi,\mu)\hat v(\xi)\tilde A(\rho,\xi)\right]d\mu.
$$
Consider the function $R_\rho(\mu):=(P_+(\mu)-P_-(\mu))\tilde d_2(\rho,\mu)\hat v(\rho)$ (where $\rho\in\check Z_k$ considered as a parameter). By virtue of lemma 4.5 one has for nonreal $\mu\in \mathbb{C}\setminus(\mathbb{R}\cup\check Z_k)$:
$$
(CR_\rho)(\mu)=P(\mu)\tilde d_2(\rho,\mu)\hat v(\rho)-
 d_2(\rho,\mu)\hat v(\rho)+\sum\limits_{\xi\in \check Z_k} d_2(\xi,\mu)\hat v(\xi)\tilde A(\rho,\xi).
$$
Taking the limit as $\mbox{Im}\mu\to +0$ we obtain for real $\mu$:
$$
(C_+R_\rho)(\mu)=P_+(\mu)\tilde d_2(\rho,\mu)\hat v(\rho)-
 d_2(\rho,\mu)\hat v(\rho)+\sum\limits_{\xi\in \check Z_k} d_2(\xi,\mu)\hat v(\xi)\tilde A(\rho,\xi).
$$
Thus, we can write:
$$
\mathcal{A}_{dd} \vp(\rho)=\frac{1}{2\pi i}\int\limits_{-\infty}^\infty f(\mu)\tilde P_+(\mu)
\left[(C_+R_\rho)(\mu)-P_+(\mu)\tilde d_2(\rho,\mu)\hat v(\rho)\right]d\mu=
$$
$$
\frac{1}{2\pi i}\int\limits_{-\infty}^\infty f(\mu)\tilde P_+(\mu)
(C_+R_\rho)(\mu)d\mu-\frac{1}{2\pi i}\int\limits_{-\infty}^\infty f(\mu)\tilde d_2(\rho,\mu)\hat v(\rho)d\mu.
$$
Using the relation
$$
\int\limits_{-\infty}^\infty  F_1(\mu) (C_+F_2)(\mu)d\mu=-\int\limits_{-\infty}^\infty  (C_-F_1)(\mu) F_2(\mu)d\mu
$$
we rewrite this as follows:
$$
(\mathcal{A}_{dd} \mathcal{B}_{dr} f)\rho=\mathcal{A}_{dd} \vp(\rho)=-\frac{1}{2\pi i}\int\limits_{-\infty}^\infty \left[C_-(f\tilde P_+)\right](\mu)
R_\rho(\mu)d\mu-\frac{1}{2\pi i}\int\limits_{-\infty}^\infty f(\mu)\tilde d_2(\rho,\mu)\hat v(\rho)d\mu=
$$
$$
-\frac{1}{2\pi i}\int\limits_{-\infty}^\infty \left[C_-(f\tilde P_+)\right](\mu)(P_+(\mu)-P_-(\mu))\tilde d_2(\rho,\mu)\hat v(\rho)
d\mu-\frac{1}{2\pi i}\int\limits_{-\infty}^\infty f(\mu)\tilde d_2(\rho,\mu)\hat v(\rho)d\mu.
$$
Now let $\vp=\mathcal{B}_{rr} f$. We rewrite it as $\vp(\mu)=f(\mu)+C_-\left[(f\tilde P_+)\right](\mu)(P_+(\mu)-P_-(\mu))$. Then
$$
(\mathcal{A}_{dr} \mathcal{B}_{rr} f)(\rho)=(\mathcal{A}_{dr} \vp)(\rho)=$$$$\frac{1}{2\pi i}\int\limits_{-\infty}^\infty f(\mu)\tilde d_2(\rho,\mu)\hat v(\rho)d\mu+
\frac{1}{2\pi i}\int\limits_{-\infty}^\infty \left[C_-(f\tilde P_+)\right](\mu)(P_+(\mu)-P_-(\mu))\tilde d_2(\rho,\mu)\hat v(\rho)
d\mu.
$$
Analogous calculations yield
$$
(\mathcal{B}_{dd} \mathcal{A}_{dr} \vp)(\rho)=\frac{1}{2\pi i}\int\limits_{-\infty}^\infty \left[C_-\vp\right](\mu)(\tilde P_+(\mu)-\tilde P_-(\mu)) d_2(\rho,\mu)\hat v(\rho)
d\mu+\frac{1}{2\pi i}\int\limits_{-\infty}^\infty \vp(\mu)\tilde P_+(\mu) d_2(\rho,\mu)\hat v(\rho)d\mu.
$$
On the other hand, if $f=\mathcal{A}_{rr} \vp$ then $(f\tilde P_+)(\mu)=\vp(\mu)\tilde P_+(\mu)+[C_-\vp](\mu)(\tilde P_+(\mu)-\tilde P_-(\mu))$ and
$$
(\mathcal{B}_{dr} \mathcal{A}_{rr} \vp)(\rho)=(\mathcal{B}_{dr} f)(\rho)=$$$$-\frac{1}{2\pi i}\int\limits_{-\infty}^\infty \vp(\mu)\tilde P_+(\mu) d_2(\rho,\mu)\hat v(\rho)d\mu
-\frac{1}{2\pi i}\int\limits_{-\infty}^\infty \left[C_-\vp\right](\mu)(\tilde P_+(\mu)-\tilde P_-(\mu)) d_2(\rho,\mu)\hat v(\rho)
d\mu.
$$
Thus, we have $\mathcal{B}_{dd} \mathcal{A}_{dr} +\mathcal{B}_{dr} \mathcal{A}_{rr} =0$.

{\bf Position (d,d).} Let $\vp=\mathcal{B}_{rd} f$. Then
$$
(\mathcal{A}_{dr} \vp)(\rho)=\frac{1}{2\pi i}\int\limits_{-\infty}^\infty \left\{\sum\limits_{\mu\in\check Z_k}f(\mu)d_1(\xi,\mu)\right\}
(P_+(\xi)-P_-(\xi))\tilde d_2(\rho,\xi)\hat v(\rho)d\xi.
$$
Using lemma 4.5 we can rewrite it as follows
$$
(\mathcal{A}_{dr} \mathcal{B}_{rd} f)(\rho)=(\mathcal{A}_{dr} \vp)(\rho)=\sum\limits_{\mu\in\check Z_k}f(\mu)\left\{\tilde A(\rho,\mu)+A(\rho,\mu)-
\sum\limits_{\xi\in\check Z_k}A(\xi,\mu)\tilde A(\rho,\xi)\right\}.
$$
On the other hand we have:
$$
(\mathcal{A}_{dd} \mathcal{B}_{dd} f)(\rho)=f(\rho)-\sum\limits_{\mu\in\check Z_k} f(\mu)A(\rho,\mu)-
\sum\limits_{\xi\in\check Z_k}\left\{f(\xi)-\sum\limits_{\mu\in\check Z_k} f(\mu)A(\xi,\mu)\right\}\tilde A(\rho,\xi)
$$
and thus we obtain $(\mathcal{A}_{dr} \mathcal{B}_{rd} f)(\rho)+(\mathcal{A}_{dd} \mathcal{B}_{dd} f)(\rho)=f(\rho)$. Symmetrical calculations yield $\mathcal{B}_{dr} \mathcal{A}_{rd} +\mathcal{B}_{dd} \mathcal{A}_{dd} =E$.$\hfil\Box$

\medskip
Now we can formulate the constructive procedure for solution of the problem $IP(k)$.

\medskip
{ \bf Procedure 4.1.} Given scattering data $J_k$ (with some fixed $k\in\{1,\dots,p\}$) and $\nu_{0j}, \sigma_j, \sigma_{j1}, \sigma_{j2}$, $j=\overline{1,p}$. \\
1. Choose $\tilde L$ with the same $\nu_{0j}, \sigma_j, \sigma_{j1}, \sigma_{j2}$, $j=\overline{1,p}$ satisfying the conditions $G$, $R_0$, $R_\infty$.\\
2. Calculate $v(\rho)$, $\rho\in \mathbb{R}\cup\check Z_k$ using (4.10) and lemma 4.2.\\
3. From given $v(\rho)$, $\rho\in \mathbb{R}\cup\check Z_k$ and $\tilde L$ find $V(x,\rho)$, $\tilde d_j(x,\rho,\mu), j=1,2$, $\tilde A(x,\rho,\mu)$.\\
4. For each fixed $x>0$ find $\Phi(x,\rho), \rho\in\mathbb{R}\cup\check Z_k$ as a (unique) solution of the linear system (4.13), (4.15).\\
5. Given $\Phi(x,\rho)$ calculate $P_1(x,\rho)$, $x>0$, $\rho\in\mathbb{C}\setminus{\mathbb{R}\cup\check Z_k}$ via (4.12).\\
6. Given $P_1(x,\rho)$ find $f_k(x,\rho)=P_{11}(x,\rho)\tilde f_k(x,\rho)+P_{11}(x,\rho)\tilde f'_k(x,\rho)$, $\psi_{kk}(x,\rho)=P_{11}(x,\rho)\tilde \psi_{kk}+P_{11}(x,\rho)\tilde \psi'_{kk}(x,\rho)$.\\
7. Calculate $q_k(x)=f''_k(x,\rho)/f_k(x,\rho)+\rho^2-\nu_{0k}x^{-2}$ (where $\rho$ for each fixed $x>0$ is arbitrary such that $f_k(x,\rho)\neq 0$).

\section{Inverse scattering on the graph}

Here we consider the following "complete" inverse scattering problem.

\medskip
{\bf Problem $IP(\Gamma)$.} Given scattering data $J$ ($=\{J_k\}_{k=1}^{p-1}$) find the potential on $\Gamma$, i.e. all the functions $q_k, k=\overline{1,p}$.

\medskip
First we note that the procedure described in previous section allows us to recover (uniquely) the potentials $q_k, k=\overline{1,p-1}$. Thus, in order to complete our solution of the problem we have to recover $q_p$. Let us consider the matching conditions for  Weyl-type solution $\psi_1(\rho)$. The following relation is a direct sequence of (3.2):
$$
\sum\limits_{j=1}^p \frac{U_{j2}(\psi_{1j}(\cdot,\rho))}{U_{j1}(\psi_{1j}(\cdot,\rho))}=0. \eqno(5.1)
$$
Since for $j=\overline{2,p}$  $\psi_{kj}(x,\rho)=\gamma_{kj}(\rho)f_j(x,\rho)$, we have:
$$
\frac{U_{j2}(\psi_{1j}(\cdot,\rho))}{U_{j1}(\psi_{1j}(\cdot,\rho))}=\frac{\sigma_{j1}+\sigma_{j2}m_j(\lambda)}{\sigma_j},
$$
where $m_j(\lambda)=b_{j2}(\rho)/b_{j1}(\rho)$ are the local Weyl functions on the rays $\mathcal{R}_j$. Thus, (5.1) can be rewritten as follows:
$$
\frac{\sigma_{p1}+\sigma_{p2}m_p(\lambda)}{\sigma_p}=-\sum\limits_{j=2}^{p-1}\frac{\sigma_{j1}+\sigma_{j2}m_j(\lambda)}{\sigma_j}-
\frac{U_{12}(\psi_{11}(\cdot,\rho))}{U_{11}(\psi_{11}(\cdot,\rho))}. \eqno(5.2)
$$
Now we note that all the Weyl functions $m_j, j=\overline{2,p-1}$ are known since corresponding potentials $q_j$ were already recovered from the solution of problems $IP(j), j=\overline{2,p-1}$. Moreover, the Weyl-type solution $\psi_{11}(x,\rho)$ can also be found from the solution of the problem $IP(1)$. Thus, all the terms in the right-hand side of (5.2) are known and we can use (5.2) to find the Weyl function $m_p$ that actually completes the solution of the problem $IP(\Gamma)$.

The following theorem summarizes our results.

\medskip
{\bf Theorem 5.1.}  Specification of the scattering data $J$ determines uniquely the potential on $\Gamma$. The functions $q_k, k=\overline{1,p}$ can be recovered by the following procedure.\\
1. For $k=\overline{1,p-1}$ using the procedure 4.1 solve the problems $IP(k)$ and recover $q_k$. Calculate $\psi_{11}(x,\rho)$ (while solving the problem $IP(1)$).\\
2. Find $m_p$ from (5.2).\\
3. Given $m_p$ recover $q_p(x), x>0$ by solving the (local) inverse spectral problem on the semi-axis $x>0$ by means of the procedure described in {\cite{Yur4}}. \\

\medskip
{\bf Acknowledgement.} This work was supported by the Russian Ministry of
Education and Science (Grant 1.1436.2014K) and Russian Fund of Basic Research (Grants 13-01-00134, 15-01-04864).


\begin{thebibliography}{3}

\bibitem{FaP} Faddeev M. and Pavlov B., Model of free electrons and
the scattering problem. Teor. Mat. Fiz. 55, no.2 (1983), 257-269
(Russian); English transl. in Theor. Math. Phys. 55 (1983), 485-492.

\bibitem{Ex} Exner P., Contact interactions on graph superlattices.
J . Phys. A: Math. Gen. 29 (1996), 87-102.

\bibitem{KoS} Kottos T. and Smilansky U., Quantum chaos on graphs.
Phys. Rev. Lett. 79 (1997), 4794-4797.

\bibitem{PoB} Pokornyi Yu.V. and Borovskikh A.V., Differential
equations on networks (geometric graphs). J. Math. Sci. (N.Y.) 119,
no.6. (2004), 691-718.

\bibitem{Kuch} Kuchment P., Quantum graphs Waves Random Media 14 (2004), S107–S128.

\bibitem{AKu} Avdonin S. and Kurasov P., Inverse problems for quantum trees. Inverse Problems and Imaging 2 (2008), 1–21.

\bibitem{Ger} Gerasimenko N.I., Inverse scattering problems on a
noncompact graph. Teoret. Mat. Fiz. 74 (1988), no. 2, 187-200;
English transl. in Theor. Math. Phys. 75 (1988),
460-470.

\bibitem{BrW} Brown B.M. and Weikard R., A Borg-Levinson theorem
for trees, Proc. R. Soc. Lond. Ser. A Math. Phys. Eng. Sci. 461,
no.2062 (2005), 3231-3243.

\bibitem{Bel1} Belishev M.I., Boundary spectral inverse problem on a
class of graphs (trees) by the BC method, Inverse Problems 20
(2004), 647-672.

\bibitem{Yur1} Yurko V.A., Inverse spectral
problems for Sturm-Liouville operators on graphs, Inverse Problems
21 (2005), 1075-1086.

\bibitem{Yur2} Yurko V.A., Inverse spectral problem for differential
operators on arbitrary compact graphs, J. Inverse Ill-Posed Probl. 18(2010), no.3.

\bibitem{Troo} Trooshin I., Marchenko V. and Mochizuki K. Inverse scattering on a graph containing circle.
Analytic methods of analysis and DEs: AMADE 2006, 237–-243, Camb.
Sci. Publ., Cambridge, 2008.

\bibitem{KuS} Kurasov P. and Sternberg
F., On the inverse scattering problem on branching graphs. J.
Phys. A, 35(2002), 101-121.

\bibitem{Kur} Kurasov P., Inverse problems for Aharonov-Bohm rings, Math. Proc.
Cambridge Philos. Soc. 148 (2010), 331--362.

\bibitem{ISI} Ignatyev M, Inverse scattering problem for Sturm--Liouville operator on one-vertex noncompact graph with a cycle, Tamkan J. of Mathematics 42, N3 (2011), 365--384.

\bibitem{Tro1} Trooshin I. and Mochizuki K. Spectral problems and scattering on noncompact
star-shaped graphs containing finite rays, J. Inverse Ill-Posed Probl. 23 (2015), no. 1, 23–-40.

\bibitem{LD} L. D. Faddeev, The inverse problem in the quantum theory of scattering, Uspekhi
Mat. Nauk 14 (1959), pp. 57--119 (in Russian); Engl. transl. in J. Math.
Phys. 4 (1963), 72--104.

\bibitem{Coz76} M. Coz and C. Coudray, The Riemann solution and the inverse quantum mechanical
problem, J. Math. Phys. 17, 888–-893 (1976).

\bibitem{Coz83} M. Coz , A Marchenko equation for complex interactions with a regular analytic continuation, J. Math. Anal. Appl. 92, 66--95 (1983).

\bibitem{KoTeS} A. Kostenko, A. Sakhnovich, and G. Teschl, Inverse eigenvalue problems for perturbed spherical
Schr¨odinger operators, Inverse Problems 26, 105013, 14pp (2010).

\bibitem{Hr12} S. Albeverio, R. Hryniv, and Ya. Mykytyuk, Scattering theory for Schr¨odinger operators
with Bessel-type potentials, J. Reine und Angew. Math. 666, 83–113 (2012).

\bibitem{KoTe} A. Kostenko and G. Teschl, Spectral asymptotics for perturbed spherical Schr¨odinger operators
and applications to quantum scattering, Comm. Math. Phys. 322, 255–275 (2013).

\bibitem{Yur3} Yurko V.A., On integral transforms connected with differential operators having singularities inside the interval, Integral Transforms and Special Functions. 5(1997), no. 3-4, 309-322.

\bibitem{Fed} Fedoseev A. E. Inverse problems for differential
equations on the half-line having a singularity in an
interior point, Tamkang J. of Math. 42(2011), no. 3, 343–-354.

\bibitem{Bea} Beals R., Deift P. and Tomei C., Direct and inverse scattering on the line.
Math. Surveys and Monographs. V.28, Amer. Math. Soc, Providence: RI, 1988.

\bibitem{YuB} Yurko, V.A. Method of Spectral Mappings in the Inverse Problem Theory. Inverse and
Ill-Posed Problems Series, Utrecht: VSP (2002).

\bibitem{Yur4} Yurko V.A., On higher-order differential operators with a singular point, Inverse Problems. 9 (1993), 495--502.



\end{thebibliography}
\end{document}